\documentclass[12pt,english]{article}
\usepackage[T1]{fontenc}
\usepackage[latin9]{inputenc}
\usepackage{amssymb}

\makeatletter
\usepackage[T1]{fontenc}
\usepackage[latin9]{inputenc}
\usepackage{color}
\usepackage{amsmath}
\usepackage{graphicx}
\usepackage{amsfonts}
\usepackage{a4wide}

\usepackage{babel}%

\usepackage{mathrsfs}
\usepackage[active]{srcltx}

\providecommand{\U}[1]{\protect\rule{.1in}{.1in}}

\newtheorem{theorem}{Theorem}

\newtheorem{corollary}[theorem]{Corollary}

\newtheorem{definition}[theorem]{Definition}

\newtheorem{proposition}[theorem]{Proposition}
\newtheorem{remark}[theorem]{Remark}

\usepackage[usenames,dvipsnames]{pstricks}
\usepackage{epsfig}
\usepackage{pst-grad} 
\usepackage{pst-plot} 

\makeatother

\usepackage{babel}
\begin{document}

\title{Maximal monotone normal cones in locally convex spaces}

\author{M.D. Voisei}

\date{{}}
\maketitle
\begin{abstract}
Equivalent conditions that make the normal cone maximal monotone are
investigated in the general settings of locally convex spaces. Some
consequences such as Bishop Phelps and sum representability results
are presented in the last part. 
\end{abstract}

\section{Preliminaries }

The aim of this paper is to characterize the subsets $C$ of a locally
convex space $(X,\tau)$ whose normal cone $N_{C}$ is a maximal monotone
operator. 

Here $(X,\tau)$ is a non-trivial (that is, $X\neq\{0\}$) real Hausdorff
separated locally convex space (LCS for short), $X^{\ast}$ is its
topological dual usually endowed with the weak-star topology denoted
by $w^{*}$, $(X^{*},w^{*})^{*}$ is identified with $X$, $\left\langle x,x^{\ast}\right\rangle :=x^{\ast}(x)=:c(x,x^{\ast})$,
for $x\in X$, $x^{\ast}\in X^{\ast}$ denotes the \emph{duality product}
or \emph{coupling }of $X\times X^{\ast}$, and $\operatorname*{Graph}T=\{(x,x^{*})\in X\times X^{*}\mid x^{*}\in T(x)\}$
stands for the \emph{graph} of $T:X\rightrightarrows X^{*}$. 

Rockafellar showed in \cite[Theorem\ A]{MR0262827} that when $X$
is a Banach space and $f:X\to\overline{\mathbb{R}}$ is proper convex
lower semicontinuous then its \emph{convex subdifferential} $\partial f:X\rightrightarrows X^{*}$,
defined by $x^{*}\in\partial f(x)$ if $f(x)$ is finite and for every
$y\in X$, $f(y)\ge f(x)+\langle y-x,x^{*}\rangle$, is maximal monotone
($\partial f\in\mathfrak{M}(X)$ for short). In particular the \emph{normal
cone} to $C$ which is given by $N_{C}=\partial\iota_{C}\in\mathfrak{M}(X)$,
whenever $C\subset X$ is closed convex. Here $\iota_{C}(x)=0$, for
$x\in C$; $\iota_{C}(x)=+\infty$, for $x\in X\setminus C$ denotes
the \emph{indicator function} of $C$. 

Therefore it is interesting to find when a normal cone is maximal
monotone outside the Banach space context. 

Our main argument stems from the explicit form of the normal cone
Fitzpatrick function (see Theorem \ref{fnc} below) and our characterization
of maximal monotone operators as representable and of type NI (see
\cite[Theorem\ 2.3]{MR2207807} or \cite[Theorem\ 3.4]{MR2453098}). 

Recall that the \emph{Fitzpatrick function} $\varphi_{T}:X\times X^{*}\rightarrow\overline{\mathbb{R}}$
of a multi-valued operator $T:X\rightrightarrows X^{*}$ is given
by (see \cite{MR1009594})
\begin{equation}
\varphi_{T}(x,x^{*}):=\sup\{\langle x-a,a^{*}\rangle+\langle a,x^{*}\rangle\mid(a,a^{*})\in\operatorname*{Graph}T\},\ (x,x^{*})\in X\times X^{*}.\label{ff}
\end{equation}

As usual, given a LCS $(E,\mu)$ and $A\subset E$ we denote by ``$\operatorname*{conv}A$''
the \emph{convex hull} of $A$, ``$\operatorname*{span}A$'' the
linear hull of $A$, ``$\operatorname*{cl}_{\mu}(A)=\overline{A}^{\mu}$''
the $\mu-$\emph{closure} of $A$, ``$\operatorname*{int}_{\mu}A$''
the $\mu-$\emph{topological interior }of $A$, ``$\operatorname*{core}A$''
the \emph{algebraic interior} of $A$. The use of the  $\mu-$notation
is not enforced when the topology $\mu$ is clearly understood. 

For $f,g:E\rightarrow\overline{\mathbb{R}}$ we set $[f\leq g]:=\{x\in E\mid f(x)\leq g(x)\}$;
the sets $[f=g]$, $[f<g]$, and $[f>g]$ being defined in a similar
manner. We write $f\ge g$ shorter for $f(z)\ge g(z)$, for every
$z\in E$. 

For a multi-function $T:X\rightrightarrows X^{*}$, $D(T)=\operatorname*{Pr}_{X}(\operatorname*{Graph}T)$,
$R(T)=\operatorname*{Pr}_{X^{*}}(\operatorname*{Graph}T)$ stand for
the domain and the range of $T$ respectively, where $\operatorname*{Pr}_{X}$,
$\operatorname*{Pr}_{X^{*}}$ denote the projections of $X\times X^{*}$
onto $X$, $X^{*}$ respectively. When no confusion can occur, $T:X\rightrightarrows X^{\ast}$
will be identified with $\operatorname*{Graph}T\subset X\times X^{*}$.

The restriction of an operator $T:X\rightrightarrows X^{*}$ to $U\subset X$
is the operator $T|_{U}:X\rightrightarrows X^{*}$ defined by $\operatorname*{Graph}(T|_{U})=\operatorname*{Graph}T\cap(U\times X^{*})$. 

The operator $T^{+}:X\rightrightarrows X^{\ast}$ whose graph is $\operatorname*{Graph}(T^{+}):=[\varphi_{T}\le c]$
describes all $(x,x^{*})\in X\times X^{*}$ that are monotonically
related (m.r. for short) to $T$, that is $(x,x^{*})\in[\varphi_{T}\le c]$
iff, for every $(a,a^{*})\in\operatorname*{Graph}T$, $\langle x-a,x^{*}-a^{*}\rangle\ge0$.

We consider the following classes of functions and operators on $(X,\tau)$
\begin{description}
\item [{$\Lambda(X)$}] the class formed by proper convex functions $f:X\rightarrow\overline{\mathbb{R}}$.
Recall that $f$ is \emph{proper} if $\operatorname*{dom}f:=\{x\in X\mid f(x)<\infty\}$
is nonempty and $f$ does not take the value $-\infty$, 
\item [{$\Gamma_{\tau}(X)$}] the class of functions $f\in\Lambda(X)$
that are $\tau$\textendash lower semi-continuous (\emph{$\tau$\textendash }lsc
for short), 
\item [{$\mathcal{M}(X)$}] the class of non-empty monotone operators $T:X\rightrightarrows X^{\ast}$
($\operatorname*{Graph}T\neq\emptyset$). Recall that $T:X\rightrightarrows X^{\ast}$
is \emph{monotone} if, for all $(x_{1},x_{1}^{*}),(x_{2},x_{2}^{*})\in\operatorname*{Graph}T$,
$$\left\langle x_{1}-x_{2},x_{1}^{\ast}-x_{2}^{\ast}\right\rangle \geq0$$
or, equivalently, $\operatorname*{Graph}T\subset\operatorname*{Graph}(T^{+})$. 
\item [{$\mathfrak{M}(X)$}] the class of maximal monotone operators $T:X\rightrightarrows X^{*}$.
The maximality is understood in the sense of graph inclusion as subsets
of $X\times X^{*}$. It is easily seen that $T\in\mathfrak{M}(X)$
iff $T=T^{+}$. 
\end{description}
To a proper function $f:(X,\tau)\rightarrow\overline{\mathbb{R}}$
we associate the following notions: 
\begin{description}
\item [{$\operatorname*{Epi}f:=\{(x,t)\in X\times\mathbb{R}\mid f(x)\leq t\}$}] is
the \emph{epigraph} of $f$,
\item [{$\operatorname*{conv}f:X\rightarrow\overline{\mathbb{R}}$,}] the
\emph{convex hull} of $f$, which is the greatest convex function
majorized by $f$, $(\operatorname*{conv}f)(x):=\inf\{t\in\mathbb{R}\mid(x,t)\in\operatorname*{conv}(\operatorname*{Epi}f)\}$
for $x\in X$,
\item [{$\operatorname*{cl}_{\tau}\operatorname*{conv}f:X\rightarrow\overline{\mathbb{R}}$,}] the
\emph{$\tau-$lsc convex hull} of $f$, which is the greatest \emph{$\tau$\textendash }lsc
convex function majorized by $f$, $(\operatorname*{cl}_{\tau}\operatorname*{conv}f)(x):=\inf\{t\in\mathbb{R}\mid(x,t)\in\operatorname*{cl}_{\tau}(\operatorname*{conv}\operatorname*{Epi}f)\}$
for $x\in X$,
\item [{$f^{\ast}:X^{\ast}\rightarrow\overline{\mathbb{R}}$}] is the \emph{convex
conjugate} of $f:X\rightarrow\overline{\mathbb{R}}$ with respect
to the dual system $(X,X^{\ast})$, $f^{\ast}(x^{\ast}):=\sup\{\left\langle x,x^{\ast}\right\rangle -f(x)\mid x\in X\}$
for $x^{\ast}\in X^{\ast}$.
\end{description}
Accordingly, $\sigma_{C}(x^{*}):=\sup\{\langle x,x^{*}\rangle\mid x\in C\}=\iota_{C}^{*}(x^{*})$,
for $x^{*}\in X^{*}$. Recall that $f^{**}:=(f^{*})^{*}=\operatorname*{cl}\operatorname*{conv}f$
whenever $\operatorname*{cl}\operatorname*{conv}f$ (or equivalently
$f^{*}$) is proper, where for functions defined in $X^{*}$, the
conjugates are taken with respect to the dual system $(X^{*},X)$. 

Throughout this article the conventions $\infty-\infty=\infty$, $\sup\emptyset=-\infty$,
and $\inf\emptyset=\infty$ are enforced while the use of the topology
notation is avoided when the topology is clearly understood.

All the considerations and results of this paper can be done with
respect to a separated dual system of vector spaces $(X,Y)$. 

\section{Support points and the maximality of the normal cone}

Given $(X,\tau)$ a LCS and $C\subset X$, we denote by 
\[
\operatorname*{Supp}C=\{x\in C\mid N_{C}(x)\neq\{0\}\}
\]
the set of \emph{support points of }$C$ and by 

\[
C^{\#}:=\{x\in X\mid\forall(a,a^{*})\in\operatorname*{Graph}N_{C},\ \langle x-a,a^{*}\rangle\le0\}
\]
the \emph{portable hull of }$C$ which is the intersection of all
the supporting half-spaces that contain $C$ and are supported at
points in $C$. 
\noindent \begin{center}
\scalebox{.6}
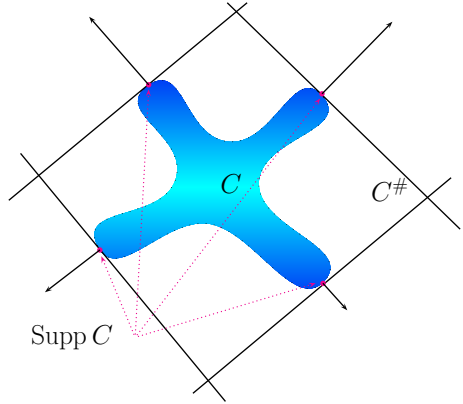
{\begin{pspicture}(0,-4.44)(10.02,4.44)
\psbezier[linewidth=0.02,linestyle=dotted,dotsep=0.16cm,fillstyle=gradient,gradlines=2000,gradmidpoint=1.5](3.08,2.58)(3.888053,3.1691098)(3.883005,1.4174647)(4.88,1.34)(5.876995,1.2625352)(6.322638,3.1105645)(6.96,2.34)(7.597362,1.5694356)(5.5612063,1.4197751)(5.54,0.42)(5.5187936,-0.5797751)(7.6558867,-0.9851405)(7.0,-1.74)(6.3441133,-2.4948595)(5.411458,-0.61042804)(4.42,-0.48)(3.4285421,-0.349572)(2.3938525,-1.860984)(1.96,-0.96)(1.5261476,-0.05901605)(3.4171054,-0.439151)(3.7,0.52)(3.9828947,1.479151)(2.2719471,1.9908901)(3.08,2.58) 

\psline[linewidth=0.03cm](5.28,4.42)(0.0,0.04) \psline[linewidth=0.03cm](0.26,1.08)(4.8,-4.42) \psline[linewidth=0.03cm](4.08,-4.24)(9.8,0.54) \psline[linewidth=0.03cm](10.0,-0.6)(4.84,4.42) 

\usefont{T1}{ptm}{m}{n} 
\rput(8.444531,0.29){\Large $C^{\#}$} 
\usefont{T1}{ptm}{m}{n} 
\rput(4.9245315,0.41){\Large $C$} 

\rput(1.4,-3){\Large $\operatorname*{Supp}C$} 
\psdots[linecolor=magenta,dotsize=0.12](6.94,2.4) 
\psdots[linecolor=magenta,dotsize=0.12](3.1,2.6) 
\psdots[linecolor=magenta,dotsize=0.12](2.02,-1.06) 
\psdots[linecolor=magenta,dotsize=0.12](6.96,-1.8)
\psline[linestyle=dotted,dotsep=0.06cm,linecolor=magenta,linewidth=0.03cm,arrowsize=0.05cm 2.0,arrowlength=1.4,arrowinset=0.4]{->}(2.8,-3)(6.9,2.3)
\psline[linestyle=dotted,dotsep=0.06cm,linecolor=magenta,linewidth=0.03cm,arrowsize=0.05cm 2.0,arrowlength=1.4,arrowinset=0.4]{->}(2.8,-3)(3.1,2.5)
\psline[linestyle=dotted,dotsep=0.06cm,linecolor=magenta,linewidth=0.03cm,arrowsize=0.05cm 2.0,arrowlength=1.4,arrowinset=0.4]{->}(2.8,-3)(2.05,-1.2)
\psline[linestyle=dotted,dotsep=0.06cm,linecolor=magenta,linewidth=0.03cm,arrowsize=0.05cm 2.0,arrowlength=1.4,arrowinset=0.4]{->}(2.8,-3)(6.8,-1.8)

\psline[linewidth=0.03cm,arrowsize=0.05291667cm 2.0,arrowlength=1.4,arrowinset=0.4]{->}(3.1,2.6)(1.78,4.11)
\psline[linewidth=0.03cm,arrowsize=0.05291667cm 2.0,arrowlength=1.4,arrowinset=0.4]{->}(6.94,2.4)(8.5,4)
\psline[linewidth=0.03cm,arrowsize=0.05291667cm 2.0,arrowlength=1.4,arrowinset=0.4]{->}(2.02,-1.06)(0.8,-2)
\psline[linewidth=0.03cm,arrowsize=0.05291667cm 2.0,arrowlength=1.4,arrowinset=0.4]{->}(6.96,-1.8)(7.5,-2.4)

\end{pspicture}} 
\par\end{center}

From their definitions, $\operatorname*{Supp}\emptyset=\operatorname*{Supp}X=\emptyset$,
$\emptyset^{\#}=X$, and 
\[
x\in C^{\#}\Leftrightarrow\forall a\in\operatorname*{Supp}C,\forall a^{*}\in N_{C}(a),\ \langle x-a,a^{*}\rangle\le0,
\]
with the remarks that, for every $C\subset X$, one has $C^{\#}\neq\emptyset$
and $\operatorname*{cl}_{\tau}\operatorname*{conv}C\subset(\operatorname*{cl}_{\tau}\operatorname*{conv}C)^{\#}\subset C^{\#}$;
while $C^{\#}=X$ when $\operatorname*{Supp}C=\emptyset$, e.g., $X^{\#}=X$. 

\begin{theorem} \label{fnc} Let $X$ be a LCS. For every $C\subset X$,
$(x,x^{*})\in X\times X^{*}$
\begin{equation}
\varphi_{N_{C}}(x,x^{*})=\iota_{C^{\#}}(x)+\sigma_{C}(x^{*}).\label{eq:fnc}
\end{equation}
In particular, for every $\emptyset\neq C\subset X$, $(x,x^{*})\in X\times X^{*}$,
$\varphi_{N_{C}}(x,0)=\iota_{C^{\#}}(x)\ge0$, $\varphi_{N_{C}}(x,x^{*})\le\iota_{C}(x)+\sigma_{C}(x^{*})$
and 
\begin{equation}
C^{\#}:=\operatorname*{Pr}\,\!\!_{X}(\operatorname*{dom}\varphi_{N_{C}})=\{x\in X\mid\varphi_{N_{C}}(x,0)=(\le)0\}=D(N_{C}^{+}).\label{eq:}
\end{equation}
\end{theorem}

\strut

The last part of Theorem \ref{fnc} says that any monotone extension
$T$ of $N_{C}$ has $D(T)\subset C^{\#}$ or that we can extend $N_{C}$
monotonically only inside $C^{\#}\times X^{*}$. 

\begin{theorem} \label{ncmm} Let $X$ be a LCS and let $C\subset X$.
The following are equivalent

\medskip

\emph{(i)} $N_{C}\in\mathfrak{M}(X)$, 

\medskip

\emph{(ii) }$\varphi_{N_{C}}(x,x^{*})=\iota_{C}(x)+\sigma_{C}(x^{*})$,
$(x,x^{*})\in X\times X^{*}$,

\medskip

\emph{(iii)} $C^{\#}\subset(=)C$

\medskip

\emph{(iv)} $C=\{x\in X\mid\varphi_{N_{C}}(x,0)\le0\}$. \end{theorem}

Concerning the previous result note that $C$ is non-empty closed
convex whenever $N_{C}\in\mathfrak{M}(X)$. Also, subpoint (iv) can
be restated equivalently as 

\medskip

(iv)$'$ $C$ is non-empty and $C=(\supset)\{x\in X\mid\varphi_{N_{C}}(x,0)\le(=)0\}$. 

\strut

Theorem \ref{ncmm} has strong ties with the separation theorem. Assume
that $C\subset X$ is closed and convex. Then for every $x\in X\setminus C$
there is $n^{*}\in X^{*}\setminus\{0\}$ such that $\langle x,n^{*}\rangle>\sigma_{C}(n^{*})=\sup\{\langle u,n^{*}\rangle\mid u\in C\}$.
In the next theorem we see that the maximality of $N_{C}$ is equivalent
to the possibility of picking, in the previous separation argument,
of a non-zero $n^{*}$ that attains its global maximum on $C$, that
is, $n^{*}\in R(N_{C})$ which is also called a \emph{support functional}
of $C$ (see e.g. \cite{MR0154092}). 

\begin{theorem} \label{ncsep} Let $X$ be a LCS and let $C\subset X$.
Then $N_{C}\in\mathfrak{M}(X)$ (or $C=C^{\#}$) iff for every $x\in X\setminus C$
there is $n^{*}\in D(\partial\sigma_{C})=R(N_{C})$ such that $\langle x,n^{*}\rangle>\sigma_{C}(n^{*})$.
\end{theorem} 

\begin{corollary} \label{Bic} Let $X$ be a LCS and let $C\subset X$
be non-empty closed and convex. If $X$ is a Banach space or $\operatorname*{int}C\neq\emptyset$,
or $C$ is weakly compact then $N_{C}\in\mathfrak{M}(X)$. \end{corollary}

Note that $C^{\#}$, the portable hull of $C$, is the smallest set
formed by intersecting half-spaces that contain $C$ and are supported
at points in $C$, and, at the same time, the largest set on which
the normal cone $N_{C}$ can be extended monotonically. 

\begin{proposition} \label{eNC} Let $X$ be a LCS. For every $C\subset X$,
$N_{C^{\#}}|_{C}=N_{C}$, $\operatorname*{Graph}N_{C}\subset\operatorname*{Graph}N_{C^{\#}}$,
$N_{C^{\#}}\in\mathfrak{M}(X)$, and $C^{\#\#}:=(C^{\#})^{\#}=C^{\#}$.
\end{proposition}

\begin{remark}\emph{ Every (maximal) monotone extension of $N_{C}$
has the domain contained in $C^{\#}$, that is, $N_{C}\subset T\in\mathcal{M}(X)\Rightarrow D(T)\subset D(N_{C}^{+})=C^{\#}$.
Therefore $N_{C^{\#}}$ is a maximal monotone extension of $N_{C}$
with the largest possible domain. In general, $N_{C^{\#}}$ is not
the only maximal monotone extension of $N_{C}$ (even with the largest
possible domain or with a normal cone structure). For example, for
$C=(0,1]\subset\mathbb{R}$, $N_{C}$ admits $N_{\overline{C}}$ and
$N_{C^{\#}}$ as two different maximal monotone extensions; moreover
$N_{C}$ has an infinity of maximal monotone extensions with the largest
possible domain $C^{\#}=(-\infty,1]$.} \end{remark}

The maximality of the subdifferential allows us to reprove and extend
some of the Bishop-Phelps results (see \cite{MR0154092}).

\begin{theorem} Let $X$ be a LCS. If, for every closed convex $C\subset X$,
$N_{C}\in\mathfrak{M}(X)$, then, for every closed convex $C\subset X$,
$\operatorname*{Supp}C$ is dense in $\operatorname*{bd}C$. \end{theorem}

\begin{remark} The previous results still holds for a fixed closed
convex $C\subset X$ if, for every $U\subset X$ closed convex such
that $C\cap\operatorname*{int}U\neq\emptyset$, $N_{C\cap U}\in\mathfrak{M}(X)$.
\end{remark}

\begin{theorem} Let $X$ be a LCS and let $C\subset X$ be closed
convex, free of lines, and finite-dimensional, that is, $\operatorname*{dim}(\operatorname*{span}C)<\infty$.
Then $N_{C}\in\mathfrak{M}(X)$ and $\operatorname*{cl}_{w^{*}}R(N_{C})=\operatorname*{cl}_{w^{*}}(\operatorname*{dom}\sigma_{C})$.

If, in addition, $C$ is bounded, then $\operatorname*{cl}_{w^{*}}R(N_{C})=X^{*}$.
\end{theorem}

\begin{proposition} under review \end{proposition}

\begin{corollary} Let $X$ be a LCS such that for every $f\in\Gamma(X)$,
$\operatorname*{Graph}\partial f\neq\emptyset$. Then for every closed
convex bounded $C\subset X$ the support functionals of $C$ are weak-star
dense in $X^{*}$. \end{corollary}

\begin{corollary} Let $(X,\|\cdot\|)$ be a Banach space. Then for
every closed convex bounded $C\subset X$ the support functionals
of $C$ are strongly dense in $X^{*}$.

If $C\subset X$ is merely closed and convex then the support functionals
of $C$ are strongly dense in the domain of $\sigma_{C}$. In other
words for every $\epsilon>0$, $f\in X^{*}$ such that $\sup_{C}f<+\infty$
there is $g\in X^{*}$, $x_{0}\in C$ such that $g(x_{0})=\sup_{C}g$
(that is, $g\in N_{C}(x_{0})$) and $\|f-g\|\le\epsilon$. \end{corollary}

\section{Partial portable hulls and sum representability}

Let $(X,\tau)$ be a LCS and $C,S\subset X$. We denote by $C_{S}^{\#}$
the \emph{partial portable hull of $C$ on $S$} which is the intersection
of all the (supporting) half-spaces that contain $C$ and are supported
at points in $S$, 
\[
C_{S}^{\#}:=\{x\in X\mid\forall a\in C\cap S,\ a^{*}\in N_{C}(a),\ \langle x-a,a^{*}\rangle\le0\}.
\]
Equivalently 
\[
x\in C_{S}^{\#}\Leftrightarrow\forall a\in S\cap\operatorname*{Supp}C,\ a^{*}\in N_{C}(a),\ \langle x-a,a^{*}\rangle\le0,
\]
with the remarks that $C_{S}^{\#}=X$ when $S\cap\operatorname*{Supp}C=\emptyset$
and that, for every $C,S\subset X$, $C_{S}^{\#}\neq\emptyset$, $C\subset\operatorname*{cl}_{\tau}\operatorname*{conv}C\subset C^{\#}\subset C_{S}^{\#}$,
$C_{S}^{\#}=C_{S\cap C}^{\#}=C_{S\cap\operatorname*{Supp}C}^{\#}$. 

Note also that $C_{S}^{\#}=C^{\#}$ whenever $S\supset\operatorname*{Supp}C$. 
\noindent \begin{center}
\scalebox{.6}
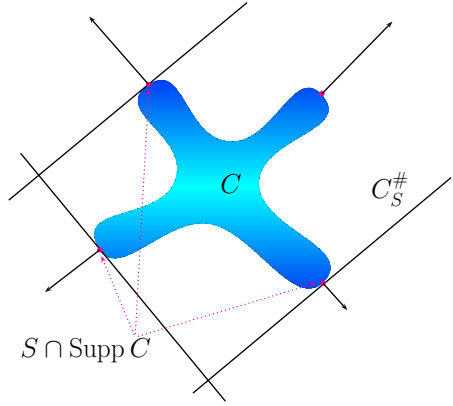
{\begin{pspicture}(0,-4.44)(10.02,4.44)
\psbezier[linewidth=0.02,linestyle=dotted,dotsep=0.16cm,fillstyle=gradient,gradlines=2000,gradmidpoint=1.5](3.08,2.58)(3.888053,3.1691098)(3.883005,1.4174647)(4.88,1.34)(5.876995,1.2625352)(6.322638,3.1105645)(6.96,2.34)(7.597362,1.5694356)(5.5612063,1.4197751)(5.54,0.42)(5.5187936,-0.5797751)(7.6558867,-0.9851405)(7.0,-1.74)(6.3441133,-2.4948595)(5.411458,-0.61042804)(4.42,-0.48)(3.4285421,-0.349572)(2.3938525,-1.860984)(1.96,-0.96)(1.5261476,-0.05901605)(3.4171054,-0.439151)(3.7,0.52)(3.9828947,1.479151)(2.2719471,1.9908901)(3.08,2.58) 

\psline[linewidth=0.03cm](5.28,4.42)(0.0,0.04) 
\psline[linewidth=0.03cm](0.26,1.08)(4.8,-4.42) 
\psline[linewidth=0.03cm](4.08,-4.24)(9.8,0.54) 

\usefont{T1}{ptm}{m}{n} 
\rput(8.444531,0.29){\Large $C^{\#}_S$} 
\usefont{T1}{ptm}{m}{n} 
\rput(4.9245315,0.41){\Large $C$} 

\rput(1.7,-3.3){\Large $S\cap \operatorname*{Supp}C$} 
\psdots[linecolor=magenta,dotsize=0.12](6.94,2.4) 
\psdots[linecolor=magenta,dotsize=0.12](3.1,2.6) 
\psdots[linecolor=magenta,dotsize=0.12](2.02,-1.06) 
\psdots[linecolor=magenta,dotsize=0.12](6.96,-1.8)
\psline[linestyle=dotted,dotsep=0.06cm,linecolor=magenta,linewidth=0.03cm,arrowsize=0.05cm 2.0,arrowlength=1.4,arrowinset=0.4]{->}(2.8,-3)(3.1,2.5)
\psline[linestyle=dotted,dotsep=0.06cm,linecolor=magenta,linewidth=0.03cm,arrowsize=0.05cm 2.0,arrowlength=1.4,arrowinset=0.4]{->}(2.8,-3)(2.05,-1.2)
\psline[linestyle=dotted,dotsep=0.06cm,linecolor=magenta,linewidth=0.03cm,arrowsize=0.05cm 2.0,arrowlength=1.4,arrowinset=0.4]{->}(2.8,-3)(6.8,-1.8)

\psline[linewidth=0.03cm,arrowsize=0.05291667cm 2.0,arrowlength=1.4,arrowinset=0.4]{->}(3.1,2.6)(1.78,4.11)
\psline[linewidth=0.03cm,arrowsize=0.05291667cm 2.0,arrowlength=1.4,arrowinset=0.4]{->}(6.94,2.4)(8.5,4)
\psline[linewidth=0.03cm,arrowsize=0.05291667cm 2.0,arrowlength=1.4,arrowinset=0.4]{->}(2.02,-1.06)(0.8,-2)
\psline[linewidth=0.03cm,arrowsize=0.05291667cm 2.0,arrowlength=1.4,arrowinset=0.4]{->}(6.96,-1.8)(7.5,-2.4)

\end{pspicture}} 
\par\end{center}

\begin{proposition} \label{ncs} Let $X$ be a LCS and let $C,S\subset X$.
Then 

\emph{(i)} $\operatorname*{Graph}(N_{C}|_{S})\subset\operatorname*{Graph}N_{C_{S}^{\#}}$
and
\[
N_{C}|_{S}=N_{C_{S}^{\#}}|_{S}\ {\rm iff}\ S\cap C=S\cap C_{S}^{\#}.
\]
In particular $N_{C^{\#}}|_{C}=N_{C}$. 

\emph{(ii)} $(C_{S}^{\#})_{S}^{\#}=(C_{S}^{\#})^{\#}=C_{S}^{\#}$.
In particular $N_{C_{S}^{\#}}\in\mathfrak{M}(X)$. \end{proposition}

\begin{theorem} under review \end{theorem}

Recall the following notion ( see Definition 14 in \cite{liro-arxiv})

\begin{definition} \emph{Let $(X,\tau)$ be a LCS. An operator }$T:X\rightrightarrows X^{*}$\emph{
is }representable\emph{ }in\emph{ $C\subset X$ or $C-$}representable\emph{
if $C\cap D(T)\ne\emptyset$ and there is $h\in\Gamma_{\tau\times w^{*}}(X\times X^{*})$
such that $h\ge c$ and $[h=c]\cap C\times X^{*}=\operatorname*{Graph}(T|_{C})$.}
\end{definition}

Using the partial portable hull we can recover the representability
of the sum between a representable operator and the normal cone (see
Theorem 35 in \cite{liro-arxiv})

\begin{theorem} \label{psi-Tnc} Let $X$ be a LCS, let $T:X\rightrightarrows X^{*}$,
and let $C\subset X$ be closed convex such that $D(T)\cap\operatorname*{int}C\neq\emptyset$.
Then $\psi_{T+N_{C}}=\psi_{T|_{C}}\square_{2}\sigma_{C_{D(T)}^{\#}\times\{0\}}$
or 
\begin{equation}
\psi_{T+N_{C}}(x,x^{*})=\min\{\psi_{T|_{C}}(x,x^{*}-u^{*})+\sigma_{C_{D(T)}^{\#}}(u^{*})\mid u^{*}\in X^{*}\},\ (x,x^{*})\in X\times X^{*},\label{psiTNC}
\end{equation}
\begin{equation}
[\psi_{T+N_{C}}=c]=[\psi_{T|_{C}}=c]+N_{C}.\label{TNCrepres}
\end{equation}
In particular, if, in addition, $T$ is $C-$representable then $T+N_{C}$
is representable. \end{theorem}

\section{Concluding remarks}

Let $(X,\tau)$ be a LCS and let us call a set $C\subset X$ \emph{portable}
if $C=C^{\#}$ (or $N_{C}\in\mathfrak{M}(X)$). Some of the results
of this paper can be summarized as follows:
\begin{itemize}
\item $C\subset X$ is portable iff for every $x\in X\setminus C$ there
is $n^{*}\in D(\partial\sigma_{C})=R(N_{C})$ such that $\langle x,n^{*}\rangle>\sigma_{C}(n^{*})$
or, equivalently, $[\sigma_{C-x}<0]\cap R(N_{C})\neq\emptyset$. 
\item If $C\subset X$ is closed convex and $X$ is a Banach space or $\operatorname*{int}C\neq\emptyset$,
or $C$ is weakly compact then $C$ is portable.
\item For every $C,S\subset X$, $C_{S}^{\#}$ is portable.
\item If every closed convex $C\subset X$ is portable then for every closed
convex $C\subset X$, $\operatorname*{Supp}C$ is dense in $\operatorname*{bd}C$. 
\item Every closed convex, free of lines, and finite-dimensional set is
portable. 
\end{itemize}
The results in this article are an incentive for studying the following
problems:
\begin{description}
\item [{(P1)}] Find cha racterizations of all LCS $X$ with the property
that, for every closed convex $C\subset X$, $N_{C}\in\mathfrak{M}(X)$. 
\item [{(P2)}] Find characterizations of all LCS $X$ with the property
that, for every $f\in\Gamma(X)$, $\operatorname*{Graph}\partial f\neq\emptyset$.
\item [{(P3)}] Find characterizations of all LCS $X$ with the property
that, for every $f\in\Gamma(X)$, $\partial f\in\mathfrak{M}(X)$.
\end{description}

\bibliographystyle{plain}

\begin{thebibliography}{1}

\bibitem{MR0154092}
Errett Bishop and R.~R. Phelps.
\newblock The support functionals of a convex set.
\newblock In {\em Proc. {S}ympos. {P}ure {M}ath., {V}ol. {VII}}, pages 27--35.
  Amer. Math. Soc., Providence, R.I., 1963.

\bibitem{MR1009594}
Simon Fitzpatrick.
\newblock Representing monotone operators by convex functions.
\newblock In {\em Workshop/{M}iniconference on {F}unctional {A}nalysis and
  {O}ptimization ({C}anberra, 1988)}, volume~20 of {\em Proc. Centre Math.
  Anal. Austral. Nat. Univ.}, pages 59--65. Austral. Nat. Univ., Canberra,
  1988.

\bibitem{MR0262827}
R.~T. Rockafellar.
\newblock On the maximal monotonicity of subdifferential mappings.
\newblock {\em Pacific J. Math.}, 33:209--216, 1970.

\bibitem{MR2207807}
M.~D. Voisei.
\newblock A maximality theorem for the sum of maximal monotone operators in
  non-reflexive {B}anach spaces.
\newblock {\em Math. Sci. Res. J.}, 10(2):36--41, 2006.

\bibitem{MR2453098}
M.~D. Voisei.
\newblock The sum and chain rules for maximal monotone operators.
\newblock {\em Set-Valued Anal.}, 16(4):461--476, 2008.

\bibitem{liro-arxiv}
M.D. Voisei.
\newblock Location, identification, and representability of monotone operators
  in locally convex spaces.
\newblock 2016.

\end{thebibliography}

\end{document}